\documentclass[reqno]{amsart}
\usepackage{amsmath,amssymb,amscd,verbatim}

\newcommand{\Z}{\ensuremath{\mathbb{Z}}}

\newcommand{\Hom}{\ensuremath{\operatorname{Hom}}}

\newcommand{\g}{\ensuremath{\mathfrak{g}}}

\newcommand{\C}{\ensuremath{\mathbb{C}} }

\newcommand{\p}[1]{\ensuremath{\bar {#1}}}

\newcommand{\EK}{\mathcal{U}_{h}^{EK}}

\newcommand{\Mb}{\overline{M}}

\newcommand{\M}{\ensuremath{\mathcal{M_{\g}}}}
\newcommand{\Mp}{\ensuremath{\mathcal{M_{\g'}}}}
\newcommand{\A}{\ensuremath{\mathcal{A}}}
\newcommand{\flip}{\sigma}

\newcommand{\dpi}{\Lambda}
\newcommand{\DJgA}{\ensuremath{U_{h}^{DJ}(\g,D) } }

\newcommand{\oneq}{1^q}

\newcommand{\gp}{g^+}
\newcommand{\gm}{g^-}
\newcommand{\Umod}{\bar{U}}

\newtheorem{Df}{Definition}
\newtheorem{theorem}[Df]{Theorem}

\newtheorem{proposition}[Df]{Proposition}

\newtheorem{lemma}[Df]{Lemma}

\newtheorem{corollary}[Df]{Corollary}
\numberwithin{Df}{section}
\numberwithin{equation}{section}

\begin{document}
\title{Monodromy of trigonometric KZ equations}
    
\author{Pavel Etingof}
\address{Department of Mathematics\\
Massachusetts Institute of Technology\\
Cambridge, Massachusetts 02139-4307}
\email{etingof@math.mit.edu}    
\author{Nathan Geer}
\address{School of Mathematics\\ 
Georgia Institute of Technology\\ 
Atlanta, GA 30332-0160}
\email{geer@math.gatech.edu}
\date{\today}

\begin{abstract}
The famous Drinfeld-Kohno theorem for simple Lie algebras states
that the monodromy representation of the Knizhnik-Zamolodchikov 
equations for these Lie algebras expresses explicitly via 
R-matrices of the corresponding Drinfeld-Jimbo quantum groups. 
This result was generalized by the second author to 
simple Lie superalgebras of type A-G. In this paper, we
generalize the Drinfeld-Kohno theorem to the case of the
trigonometric Knizhnik-Zamolodchikov equations for simple Lie
superalgebras of type A-G. The equations contain a classical
r-matrix on the Lie superalgebra, and the answer expresses 
through the quantum R-matrix of the corresponding quantum group. 

The proof is based on the quantization theory for Lie bialgebras
developed by the first author and D. Kazhdan.    
\end{abstract}

\maketitle

\section{Introduction}
Let $A$ be an associative superalgebra (i.e. an 
algebra with a $\mathbb{Z}_2$-grading), and $r\in A\otimes A$ be an element 
satisfying the (superversion of) the classical Yang-Baxter
equation:
$$
[r_{12},r_{13}]+[r_{12},r_{23}]+[r_{13},r_{23}]=0.
$$
Let ${\mathbf r}(u)=\frac{r_{21}e^u+r}{e^u-1}$.
Then $ {\mathbf r}(u)$ is a unitary classical r-matrix with 
a spectral parameter, i.e. $ {\mathbf r}(u)= -{\mathbf r}^{21}(-u)$, and 
$ {\mathbf r}(u)$ satisfies the classical Yang-Baxter equation 
with a spectral parameter, 
 \begin{equation}
\label{E:CYBEsp}
[{\mathbf r}_{12}(z_{1}-z_{2}),{\mathbf r}_{13}(z_{1}-z_{3})]+
[{\mathbf r}_{12}(z_{1}-z_{2}),{\mathbf r}_{23}(z_{2}-z_{3})]+
[{\mathbf r}_{13}(z_{1}-z_{3}),{\mathbf r}_{23}(z_{2}-z_{3})]=0.  
\end{equation}
Therefore, for any even element 
$s\in A$ such that $[s\otimes 1+1\otimes s,r]=0$ the system of differential 
equations 
\begin{equation}
\frac{\partial F}{\partial u_k}=\frac{h}{2\pi i}
\left(\sum_{j\ne k} {\mathbf r}_{kj}(u_k-u_j)+s^{(k)}\right)F
\end{equation}
for a function $F(u_1,...,u_n)$ with values in $A^{\otimes n}$ is consistent,
i.e. represents a flat connection in the trivial bundle.
This system is called the \emph{trigonometric 
Knizhnik-Zamolodchikov equations}, since it is analogous to the 
KZ system and contains the trigonometric r-matrix $ {\mathbf r}$.  

The main result of this paper is a computation of 
the monodromy representation of the trigonometric KZ equations,
and an explicit presentation of the answer for simple Lie superalgebras
of type A-G. 

In more detail, recall (see e.g. \cite{ES})
that the r-matrix $r$ gives rise to a Lie
superbialgebra $\g_+$ whose Drinfeld double $\g$ maps to $A$. 
The Lie bialgebra $\g$ admits a canonical quantization 
$U_h(\g)$ described in \cite{EK1,EK2}, which is
a quasitriangular quantized universal enveloping algebra. 
In particular, it has a universal R-matrix $\tilde R$.
The monodromy representation of the trigonometric KZ equation 
is described in terms of $\tilde R$. 

For Lie superalgebras of type A-G and standard r-matrix $r$ 
(for some polarization), the Hopf algebra $U_h(\g)$
and the R-matrix $\tilde R$ are known explicitly, since 
by the results of \cite{G04A,G04C}, they coincide with the Drinfeld-Jimbo  
quantum group and R-matrix. This allows us to give a completely
explicit description of the monodromy representation in this
case. This description is a trigonometric analog of the 
Drinfeld-Kohno theorem, describing the monodromy of the usual
(rational) KZ equations, which was proved for finite dimensional 
Lie superalgebras in \cite{G04C}. 

The paper is organized as follows. In Section 2, we introduce
the main definitions and notation. In Section 3, we formulate the
main results of the paper. In Section 4, we explain the
connection between the trigonometric KZ equations and the rational KZ
equations for an extended Lie superalgebra $\g'$. This connection
helps us, in Section 5, to calculate the monodromy of the
trigonometric KZ equations, and thus prove the main results
of the paper. 

\subsection*{Acknowledgments}  
The main result of this paper (for ordinary, rather than super, algebras) 
first appeared in an appendix to the first edition 
of the book by P. Etingof and O. Schiffmann ``Lectures on quantum
groups''. However, the proof contained mistakes, as was pointed
out by A. Haviv. We are very grateful to A. Haviv for many
detailed comments. The work of P.E. was  partially supported 
by the NSF grant DMS-0504847.

%  Section Preliminaries 

\section{Preliminaries}\label{S:Pre}
In this section we recall facts and definitions related to superspaces and Lie super(bi)algebras, for more details see \cite{K}.

A \emph{superspace} is a $\Z_{2}$-graded vector space $V=V_{\p 0}\oplus V_{\p 1}$ over $\C$.  We denote the parity of a homogeneous element $x\in V$ by $\p x\in \Z_{2}$.  We say $x$ is even (odd) if $x\in V_{\p 0}$ (resp. $x\in V_{\p 1}$).  Let $V$ and $W$ be superspaces.    Let $\flip_{V,W}:V\otimes W \rightarrow W \otimes V $ be the linear map given by
\begin{equation}
\label{E:Flip}
 \flip_{V,W}(v\otimes w)=(-1)^{\p v \p w}w\otimes v
\end{equation}
for homogeneous $v \in V$ and $w \in W$.  When it is clear we will write $\flip$ for $\flip_{V,W}$. 
   
   A linear morphism can be defined on homogeneous elements and then extended by linearity.  When it is clear and appropriate we will assume elements are homogeneous.    Throughout, all modules will be $\Z_{2}$-graded modules, i.e. module structures which preserve the $\Z_{2}$-grading (see \cite{K}).

A \emph{superalgebra} is a superspace $A$ with a multiplication which preserves the $\Z_{2}$-grading.  Let $A$ and $B$ be superalgebras with $1$.  The tensor product $A\otimes B$ is a superalgebra where the multiplication is given by $(a\otimes b)(a' \otimes b'):= (-1)^{\p b \p a'}a a' \otimes b b'$.  Let $a\in A$ and $x=\sum a_{i}\otimes a'_{i}\in A\otimes A$.   Denote the element $1 \otimes...\otimes 1 \otimes a\otimes 1\otimes...\otimes 1 \in A^{\otimes n}$ where $a$ is in the $j$th position by $a^{(j)}$.  Similarly, let  $x_{jk}$ be the element of $A^{\otimes n}$ given by $\sum 1\otimes...\otimes 1 \otimes a_{i} \otimes 1 \otimes...\otimes 1\otimes a'_{i}\otimes 1 \otimes ...\otimes 1$ if $j<k$ and $\sum (-1)^{\p{a}_{i} \p{a}'_{i}}1\otimes...\otimes 1 \otimes a'_{i} \otimes 1 \otimes...\otimes 1\otimes a_{i}\otimes 1 \otimes ...\otimes 1$ if $j>k$ where $a_{i}$ is in the $j$th position and $a'_{i}$ is in the $k$th position.  

A \emph{Lie superalgebra} is a superspace $\g=\g_{\p 0} \oplus \g_{\p 1}$ with a superbracket $[\: , ] :\g^{\otimes 2}~\rightarrow~\g$ that preserves the  $\Z_{2}$-grading, is super-antisymmetric ($[x,y]=-(-1)^{\p x \p y}[y,x]$), and satisfies the super-Jacobi identity (see \cite{K}).  
A \emph{Lie superbialgebra} is a Lie superalgebra $\g$ with a linear map $\delta : \g \rightarrow \wedge^{2}\g$ that preserves the $\Z_{2}$-grading and satisfies both the super-coJacobi identity and cocycle condition (see \cite{A}).  For an element $r=\sum_{m}x_m \otimes y_m\in \g \otimes \g$ let $\rho_r\in U(\g)$ be defined by $\rho_r=\sum_m x_m y_m$.

A triple $(\g,\g_{+},\g_{-})$ of finite dimensional Lie superalgebras is a finite dimensional \emph{super Manin triple} if $\g$ has a non-degenerate super-symmetric invariant bilinear form $<,>$, such that $\g \cong \g_{+}\oplus \g_{-}$ as superspaces, and $\g_{+}$ and $\g_{-}$ are isotropic Lie sub-superalgebras of $\g$.  % i.e. $<x,y>=0$ for all $x,y\in \g_{\pm}$
There is a one-to-one correspondence between finite dimensional super Manin triples and finite dimensional Lie superbialgebras (see \cite[Proposition 1]{A}). 

Let $(\g,\g_{+},\g_{-})$ be a  finite dimensional \emph{super Manin triple}.  Let $\gp_{1},...,\gp_{n}$ be a homogeneous basis of $\g_{+}$.  Let $\gm_{1},...,\gm_{n}$ be the basis of $\g_{-}$ which is dual to  $\gp_{1},...,\gp_{n}$, i.e. $<\gm_{i},\gp_{j}>=\delta_{i,j}$.  Define $r^{\g}= \sum \gp_{i}\otimes  \gm_{i} \in \g_{+}\otimes \g_{-}\subset \g \otimes \g$.  Then $r^{\g}$ satisfies \emph{classical Yang-Baxter  equation}, i.e.
$$[r^{\g}_{12},r^{\g}_{13}]+[r^{\g}_{12},r^{\g}_{23}]+[r^{\g}_{13},r^{\g}_{23}]=0.$$ 
  We call $\g$ the \emph{double} of $\g_{+}$ and denote it by $D(\g_{+})$.

% SECTION    The trigonometric KZ equations and the main results

\section{The trigonometric KZ equations and the main
results}\label{S:tKZ}

\subsection{Trigonometric KZ equations, configuration spaces and braid groups} 

Let $A$ be an associative superalgebra over $\C$ with $1$, 
and $r$ an even element of  $A\otimes A$ which satisfies
classical Yang-Baxter equation.    
Then it is well known and 
not difficult to show (see e.g. \cite{ES}) 
that the function ${\mathbf r}(u):=\frac{r_{21}e^u+r}{e^u-1}$ is a unitary classical r-matrix with 
a spectral parameter.  Therefore, for any even element 
$s\in A$ such that $[s\otimes 1+1\otimes s,r]=0$ the system of differential 
equations 
\begin{equation}\label{E:AA1}
\frac{\partial F}{\partial u_k}=\frac{h}{2\pi i}
\left(\sum_{j\ne k}{\mathbf r}_{kj}(u_k-u_j)+s^{(k)}\right)F
\end{equation}
for a function $F(u_1,...,u_n)$ with values in $A^{\otimes n}$ is consistent,
i.e. represents a flat connection in the trivial bundle.
This system is called the \emph{trigonometric Knizhnik-Zamolodchikov equations}, since it is analogous to the KZ system and contains the trigonometric r-matrix ${\mathbf r}$.  
As explained in the introduction, our goal in this paper is to compute the monodromy representation 
of this equation. 

Equation~(\ref{E:AA1}) is regular in the region $X_n^*\subset (\C/2\pi i\Z)^n$ 
defined by the inequalities $u_k\neq u_j$, $k\neq j$. 
It is also invariant under the simultaneous action 
of the symmetric group $S_n$ on the variables $u_1,...,u_n$ and 
the space of values $A^{\otimes n}$. Thus, the monodromy 
representation for~(\ref{E:AA1}) is a group homomorphism 
$\theta: \pi_1(X_n^*/S_n,x_0)\to (A^{\otimes n}[[h]])^\times$, 
where $x_0\in X_n^*/S_n$ is a fixed 
point, and $(A^{\otimes n}[[h]])^\times$ is the multiplicative group of the algebra 
$A^{\otimes n}[[h]]$. 

The representation $\theta$ is defined by the following property. Let $F_{x_0}(u)$ be the unique solution of equation~(\ref{E:AA1}) such that $F_{x_0}(x_0)=1$. Then for any path $\gamma\in \pi_1(X_n^*/S_n,x_0)$, one has $A_\gamma F(u)=F(u)\theta(\gamma)$ near $x_0$, where $A_\gamma$ denotes the analytic continuation along $\gamma$.

{\bf Remark.} Our convention for the composition law in the
fundamental group is as follows: if $g_1,g_2$ are loops, then
$g_1g_2$ is the loop obtained by first tracing out $g_2$ and then
$g_1$. 

 Let us describe the group $\Pi_n=\pi_1(X_n^*/S_n,x_0)$ 
as a braid group. For any surface $\Sigma$, 
we can define the braid group $B_n(\Sigma)$ of this surface, to be the 
fundamental group $\pi_1(X_n(\Sigma)/S_n,x_0)$, where $X_n(\Sigma)=
\{(u_1,...,u_n)\in \Sigma^n: u_k\neq u_j,k\neq j\}$ is the configuration 
space of $\Sigma$. In particular, $X_n^*=X_n(\C/2\pi i\Z)$, 
so $\Pi_n=B_n(\C/2\pi i\Z)$. 

It is convenient to identify $\C/2\pi i\Z$ with $\C^*$, 
using the change of variable $z=e^u$. Then $X_n^*$ is identified 
with the open subset in $\C^n$ defined by the equations $z_i\neq 0$, 
$z_i\neq z_j$. Therefore, the group $\Pi_n$ is isomorphic to 
the group of all elements of $B_{n+1}$ such that the end of the first braid 
(corresponding to the point $0$) coincides with its beginning. 

It is convenient to describe $\Pi_n$ by generators and relations.
Choose the initial point on $X_n(\C^*)/S_n$ to be $(1,2,...,n)$.  
First of all, we have an embedding $B_n\to \Pi_n$, induced 
by the embedding $\C^*\to \C$. Thus, $\Pi_n$ has natural generators
$b_1,...,b_{n-1}$ (namely, $b_j$ exchanges $z_j$ with $z_{j+1}$,
with $z_{j+1}$ passing above $z_j$) with the usual braid relations. Besides, we have generators
$X_j$, $j=1,...,n$, which correspond to the point $j$ going around the circle 
$|z|=j$ counterclockwise. These generators satisfy the relations
\begin{equation}\label{E:AA2}
X_jX_k=X_kX_j, b_iX_i=X_{i+1}b_i^{-1}.
\end{equation}

The following proposition is standard. 
\begin{proposition} The group $\Pi_n$ is generated by $B_n$ and $\{X_j\}$, with defining relations~(\ref{E:AA2}).
\end{proposition}

\subsection{Quantization theory and the main results}

We will use the theory of quantization of Lie bialgebras
developed by Etingof and Kazhdan in \cite{EK1,EK2} to compute the
monodromy representation $\theta$. So let us recall some basic
facts from this theory.

Let $\EK(\g_{+})$ be the E-K
quantization associated to a Lie bialgebra $\g_{+}$.  If
$\g=D(\g_{+})$ is the double of a finite dimensional Lie
bialgebra $\g_{+}$, then $\EK(\g)$ is naturally isomorphic to the
quantum double of $\EK(\g_+)$, and therefore has a  quantum R-matrix
$\tilde{R}$.  By construction, $\EK(\g)$ is canonically
isomorphic as an algebra to $U(\g)[[h]]$, and this isomorphism 
$\zeta: U(\g)[[h]]\to \EK(\g)$ is
the identity modulo $h$; it becomes an isomorphism of Hopf
algebras if the coproduct in $U(\g)[[h]]$ is twisted by a certain
pseudotwist $J$ constructed in \cite{EK1} from the Drinfeld
associator. Thus, we can assume 
that $\tilde{R}\in (U(\g)\otimes U(\g))[[h]]$.
 
Let $A$ be an associative algebra, and $r\in A\otimes A$ be a
classical r-matrix. Then from Theorem 5.1 of \cite{EK1}, 
we have that there exists a quantum R-matrix $R\in (A\otimes
A)[[h]]$ 
such that $R=1 + hr \mod h^{2}$.   
This R-matrix is constructed as follows.  
The r-matrix $r$ gives rise to a Manin triple
$(\g,\g_+,\g_-)$, and a homomorphism of associative algebras 
$\pi: U(\g)\rightarrow A$ (see Section 5 of \cite{EK1}).  
Then we set $R=(\pi \otimes \pi)(\tilde{R})$.

In \cite{G04A} the theory of the E-K quantization is generalized to the setting of Lie superbialgebras.  From this work it follows that an analogous theorem to Theorem~5.1 of \cite{EK1}  holds for Lie superbialgebras.   Therefore, for each associative superalgebra $A$ there exists a R-matrix $R\in (A\otimes A)[[h]]$ arising as in the previous paragraph.  

Our main result is the following theorem, to be proved in
the following sections. 
Let us identify $\EK(\g)$ with $U(\g)[[h]]$ using the
canonical isomorphism $\zeta$ between them. 
Let $T=(1\otimes S)(\tilde{R})$, where $S$ is the antipode of
$\EK(\g)$.  Let $C\in A^{\otimes n}[[h]]$ be the element $C:=
\pi^{\otimes n}(m_{01}(T_{01}...T_{0n}))$ where 
$m$ is the multiplication of $\EK(\g)$.  

%Main Theorem

\begin{theorem}\label{T:Main}
The monodromy representation $\theta$ of equation \eqref{E:AA1} is
isomorphic to the representation defined by the formulas 
\begin{equation}
b_i\to \flip_{ii+1}R_{ii+1}, X_1\to
e^{h(s+\rho_r)^{(1)}}C
\end{equation}
where $\rho_{r}$ is the element of $U(\g)$ defined in Section \ref{S:Pre}.
\end{theorem}

\begin{corollary}\label{elementU}
Let $U=m(T)$ be the Drinfeld element of $\EK(\g)$.
\footnote{We note that usually the term ``Drinfeld element'' is
applied to the element $S(U)$.}
Then $U=\zeta(e^{-h\rho_r})$. 
\end{corollary}

\begin{proof}
Let us apply the theorem in the case $n=1$, $s=0$, $A=U(\g)$. In this case 
the right hand side of the trigonometric KZ equation is zero, 
so the monodromy should be trivial. So from the theorem we get 
$e^{h\rho_r}m(T)=1$, which implies the statement.
\end{proof}

{\bf Remark.} 
According to Drinfeld, in a quasitriangular Hopf algebra
we have $S^2(x)=U^{-1}xU$. Thus, we have shown that in $\EK(\g)$,
one has $S^2(x)=Ad(e^{h\rho_r})$. This gives another proof of 
Proposition A3 in the Appendix to \cite{EK3}.  

We will now use Theorem \ref{T:Main} to prove a Drinfeld-Kohno like theorem for the monodromy of the trigonometric KZ equations arising from a Lie superalgebra of type A-G.   First, we will recall some important facts about Lie superalgebras.

Any two Borel subalgebras of a semisimple Lie algebra are
conjugate.  Moreover, semisimple Lie algebras are determined by
their root systems or equivalently their Dynkin diagrams. To the
contrary, not all Borel subalgebras of classical Lie
superalgebras are conjugate.  As shown by Kac \cite{K} a Lie
superalgebra can have more than one Dynkin type diagram depending
on the choice of a Borel subalgebra.  However, using Dynkin type diagrams, Kac gave a characterization of Lie superalgebras of type A-G.  

A \emph{Lie superalgebra of type A-G} is a pair $(\g,D)$ where
$\g$ is a Lie superalgebra determined by the Dynkin type diagram
$D$.   Each $D$ defines a Lie superbialgebra structure on $\g$
and an element $\bar{r}_{g}\in \g\otimes \g$ which satisfies the
classical Yang-Baxter equation (see \cite{G04C}).  Also for each
pair $(\g,D)$ Yamane \cite{Yam94} showed that there exists a
Drinfeld-Jimbo type Hopf superalgebra $\DJgA$.  
This superalgebra is given by generators and relations and has an explicit R-matrix $\bar{R}$.  

%   COROLLARY D-J for type A-G

\begin{corollary}\label{C:AG}
Let $(\g,D)$ be a Lie superalgebra of type A-G.  Let $A=U(\g)$
and $r=\bar{r}_{\g}$. Let $s$ belong to the Cartan subalgebra of
$\g$. Then the monodromy representation of
equation \eqref{E:AA1} is isomorphic to a representation of
$\Pi_n$ in $(U_h^{DJ}(\g,D)^{\otimes n})^\times$ 
\footnote{All tensor products of $k[[h]]$-modules 
are completed with respect to
the h-adic topology}
given by the formula 
\begin{equation}
b_i\to \flip_{ii+1}\bar{R}_{ii+1}, X_1\to
(e^{hs}U^{-1})^{(1)}C
\end{equation}
where $\bar{R}$ is defined above, $U=m(T)$ is the Drinfeld
element, and $C=m_{01}(T_{01}...T_{0n})$.
\end{corollary}

\begin{proof}
From \cite{EK6,G04A,G04C} we have that $\DJgA$ is isomorphic to
$\EK(\g)$ as a quasitriangular Hopf superalgebra, and this
isomorphism acts as the identity on the Cartan subalgebra. 
We also see from Corollary \ref{elementU} that under this isomorphism, 
$e^{-h\rho_r}$ goes to the Drinfeld element $U$ of
$U_h^{DJ}(\g)$. 
The corollary then follows from Theorem~\ref{T:Main}.
\end{proof}

%
% SECTION    Trigonometric and rational KZ equations

\section[Trigonometric and rational KZ equations]{The connection of the trigonometric KZ equations
with the usual (rational) KZ equations.}
Let $A$ be an associative superalgebra over $\C$ with $1$, 
and $r$ an even element of  $A\otimes A$ which satisfies classical Yang-Baxter equation.  Let $s\in A$ be an even element  such that $[s\otimes 1+1\otimes s,r]=0$.  As described above, the element $r$ gives rise to  a finite dimensional super Manin triple  $(\g,\g_{+},\g_{-})$ and a homomorphism of associative superalgebras $\pi: U(\g)\rightarrow A$.

Set $t=s+\rho_{r}\in A$.  Then $[t\otimes 1+1\otimes t,r]=0$ and
so $[t,\g_{\pm}]\subset \g_{\pm}$, i.e. $t$ induces a derivation
of the superbialgebras $\g_{\pm}$, which preserves the duality
pairing between them. Consider the Lie superbialgebra
$\g_{-}'=\g_{-}\oplus \C t$ (a semidirect product, where $t$ acts
by the corresponding derivation and $\delta(t)=0$).  Let
$\g_{+}'=\g_{+}\oplus \C t^{*}$ be the dual algebra, where
$t^{*}$ is central (even in the double).  Consider the super
Manin triple $(\g',\g'_{+},\g'_{-})$ and in particular the
Drinfeld double $\g'=D(\g_{+}')$. 

The Lie algebra and coalgebra structures on $\g'$ are as follows.
The Lie algebra structure: $t^*$ is central; commutator with $t$ is the
above derivation; commutators inside $\g_+,\g_-$ are the same as
in $\g$, and for $x\in \g_+,y\in \g_-$, 
$$
[x,y]_{\g'}=[x,y]_\g+([t,x],y)t^*. 
$$
The Lie coalgebra structure: $\delta(t)=\delta(t^*)=0$, $\delta|_{\g_-}$ is
the same as in $\g$; if $x\in \g_+$ then
$$
\delta(x)=\delta_\g(x)+[t,x]\wedge t^*.
$$ 

Let $\Umod$ be the associative superalgebra $U(\g')/(t^{*})$.
The multiplication of $U(\g')$ induces a natural  $\g'$-module
structure on $\Umod$.   Also, the homomorphism $\pi$ extends
naturally to a superalgebra homomorphism $\Umod \rightarrow A$,
where $t\mapsto s+\rho_r$ (the reason for this assignment will
become clear from the proposition below).
 Thus, to prove Theorem \ref{T:Main}, we may assume $A=\Umod$, 
$r=r_{\g}$, and $s=t-\rho_r$,
where $(\g,\g_{+},\g_{-})$ is a finite dimensional super Manin
triple, and $\Umod$ is defined as above.  

We will now define a $\g'$-module $\Mb_{-}$ as follows:  Let
$c_{(-,+1)}$ be the one dimensional $\mathbb{C} t\oplus \g_{-}\oplus \C
t^{*}$-module where $t^{*}$ acts by $1$ and $\mathbb Ct\oplus \g_{-}$ acts by 0.
Induce the  module $c_{(-,+1)}$ to the algebra $\g'$, and call this
$\g'$-module $\Mb_{-}$.  The $\g'$-module $\Mb_{+}$ is
obtained in a similar way, with $+$ and $-$ interchanged and defining the action
of $t^{*}$ to be $-1$.  From the Poincare-Birkhoff-Witt theorem
it follows that the modules $\Mb_{\pm}$ are freely generated over
$U(\g_\mp)$ by the vectors $1_\pm$ such that $\g_\pm1_\pm=0$,
$t1_\pm=0$ and $t^*1_\pm=\mp 1_\pm$.

By Frobenius reciprocity, the space $\mathrm{Hom}_{\g_{+}\oplus \g_{-}\oplus \C
t^{*}}(\Mb_-,\Mb_+^*\hat\otimes \Umod^{\otimes n})$ 
(where $\hat\otimes$ is the completed tensor product) 
is naturally isomorphic to $\Umod^{\otimes n}$.   For a function $\Psi(z_0,z_1,...,z_n)$ with values in 
$\mathrm{Hom}_{\g_{+}\oplus \g_{-}\oplus \C
t^{*}}(\Mb_-,\Mb_+^*\hat\otimes \Umod^{\otimes n})$, let
$$F_\Psi=(1_+\otimes 1^{\otimes n})(\Psi 1_-)\in \Umod^{\otimes n}.$$ 
Let $\Omega'=\Omega+t\otimes t^{*}+ t^{*}\otimes t$, where
$\Omega=r+r^{op}$, be the Casimir tensor of $\g'$. 

%  PROP  trig KZ and usual KZ

\begin{proposition}\label{P:1} $F_\Psi(0,z_1,..., z_n)$ 
satisfies the trigonometric KZ equations~(\ref{E:AA1}) for $s=t-\rho_r$
(with $z_j=e^{u_j}$)
if and only if $\Psi(0,z_1,...,z_n)$ satisfies the usual 
KZ equations with $n+1$ variables: 
\begin{equation}\label{E:AA3}
\frac{\partial \Psi}{\partial z_k}=\frac{h}{2\pi i}
\sum_{j\neq k,0\le j\leq n}\frac{\Omega'_{kj}}{z_k-z_j}\Psi \qquad (z_0=0).
\end{equation}
\end{proposition}
\begin{proof} Suppose $\Psi$ satisfies equation \eqref{E:AA3}. Then we get
\begin{equation}\label{E:AA4}
\begin{split}
\frac{2\pi i}{h}\frac{\partial F_\Psi}{\partial z_k}&=\sum_{j\neq
k}(1_+\otimes 1^{\otimes n})
\big(\frac{\Omega'_{kj}\Psi 1_-}{z_k-z_j}\big)\\
&=\sum_{j\neq k,0}\frac{\Omega'_{kj}F_\Psi}{z_k-z_j}+
(1_+\otimes 1^{\otimes n})\big(
\frac{\Omega'_{k0}\Psi 1_-}{z_k}\big).
\end{split}
\end{equation}

Since $t^{*}$ acts by $0$ on $\Umod$ then for $k\neq 0$ we have $\Omega_{kj}' F_{\Psi}=\Omega_{kj} F_{\Psi}$.   Thus, from~(\ref{E:AA4}) we get:
\begin{equation*}\label{E:AA4.5}
\begin{split}
\frac{2\pi i}{h}\frac{\partial F_\Psi}{\partial z_k}&=
\sum_{j\neq k,0}\frac{\Omega_{kj}F_\Psi}{z_k-z_j}+
(1_+\otimes 1^{\otimes n})(\frac{\Omega_{k0}'\Psi 1_-}{z_k})\\
&=
\sum_{j\neq k,0}\frac{\Omega_{kj}F_\Psi}{z_k-z_j}+(1_+\otimes 1^{\otimes n})\bigg(\frac{(t^{*(k)}t^{(0)}+t^{(k)}t^{*(0)})\Psi 1_-}{z_k} +\sum_m\frac{g_m^{+(k)}g_m^{-(0)}\Psi 1_-}{z_k}\bigg)\\
&=
\sum_{j\neq k,0}\frac{\Omega_{kj}F_\Psi}{z_k-z_j}+(1_+\otimes 1^{\otimes n})\bigg(\frac{t^{(k)}\Psi 1_-}{z_k} +\sum_m\frac{g_m^{+(k)}g_m^{-(0)}\Psi 1_-}{z_k}\bigg)
\end{split}
\end{equation*}
Now the $\g'$-invariance of $\Psi$ implies,
\begin{equation}\label{E:AA5}
\begin{split}
\frac{2\pi i}{h}\frac{\partial F_\Psi}{\partial z_k}&=
\sum_{j\neq k,0}\frac{\Omega_{kj}F_\Psi}{z_k-z_j}
+\frac{t^{(k)}F_{\Psi}}{z_{k}}-(1_+\otimes 1^{\otimes n})\bigg(
\sum_{m,1\leq j\leq n}\frac{g_m^{+(k)}g_m^{-(j)}\Psi 1_-}{z_k}\bigg).
\end{split}
\end{equation}
From~(\ref{E:AA5}), we get
\begin{equation}\label{E:AA6}
\begin{split}
\frac{2\pi i}{h}\frac{\partial F_\Psi}{\partial z_k}&=
\sum_{j\ne k,0}\frac{\Omega_{kj}F_\Psi}{z_k-z_j}
+\frac{t^{(k)}F_{\Psi}}{z_{k}}-\frac{\rho_r^{(k)}F_\Psi}{z_k}-
\sum_{m,j\neq k,j\geq 1}\frac{r_{kj}F_\Psi}{z_k}\\
&=z_k^{-1}\bigg(\sum_{j\neq k,0}{\mathbf r}^{kj}(\mathrm{ln}(z_k/z_j))+t^{(k)}-\rho_r^{(k)}\bigg)F_\Psi.
\end{split}
\end{equation}
Making the change of variable $z_i=e^{u_i}$, we reduce~(\ref{E:AA6}) to 
the trigonometric KZ equations with $s=t-\rho_r$. 
\end{proof}

%  SECTION    Monodromy of the trigonometric KZ equations

\section{Monodromy of the trigonometric KZ equations}\label{S:MTKZ}  
 Proposition \ref{P:1} allows us to use the usual KZ equations to
compute the monodromy of the trigonometric KZ equations.   Thus,
it remains to compute the monodromy of equation~(\ref{E:AA3}).  
 
Let $M_{\pm}$ be the induced modules
$\mathrm{Ind}_{\g'_{\pm}}^{\g'}c_{\pm}$, where $c_{\pm}$ is the
trivial $\g'_{\pm}$-module.  Let $\A$ be the symmetric tensor
category of topologically free super $k[[h]]$-modules, with the canonical associativity isomorphism and super commutativity isomorphism   $\flip_{V,W}$  given in \eqref{E:Flip}.
  
Following Drinfeld, let $\Mp$ be the category whose objects are $\g'$-modules and whose morphisms are given by $\Hom_{\Mp}(V,W)=\Hom_{\g'}(V,W)[[h]]$.  For any $V, W \in \Mp$, let $V \otimes W$ be the usual super tensor product.  
Let $\beta_{V,W} : V\otimes W \rightarrow W\otimes V$ be the morphism given by the action of $e^{{h \Omega'}/2}$ on $V \otimes W$ composed with the morphism $\flip_{V,W}$ given in \eqref{E:Flip}. 
Let $\Phi$ be the KZ associator arising from the usual KZ equations (see for example \cite{EK1,G04A,Kas}).  For $V,W,U\in \M$, let $\Phi_{V,W,U}$ be the morphism defined by the action of $\Phi$ on $V \otimes W\otimes U$ regarded as an element of $\Hom_{\Mp}((V \otimes W)\otimes U,V \otimes (W\otimes U))$.  The morphisms $\Phi_{V,W,U}$ and $\beta_{V,W}$ define a braided tensor structure on the category $\Mp$ (see \cite[Prop. XIII.1.4]{Kas}), which we call the \emph{Drinfeld category}.

Define $F: \Mp \rightarrow \A$ to be the functor 
$$F(V)=\Hom_{\Mp}(M_{-},M_{+}^{*}\hat\otimes V).$$
In \cite{EK1,G04A} it is shown that the 
functor $F$ is a tensor functor and that $F(V)$ is isomorphic to
$V[[h]]$ as a $k[[h]]$-module.

Let $M_{\pm }^{q}$ be the induced module
$\mathrm{Ind}_{\EK(\g'_{\pm})}^{\EK(\g')}c_{\pm}$, where $c_\pm$
is the trivial module.  

Note that by functoriality of quantization, 
we have a natural algebra isomorphism 
\begin{equation}\label{E:9}
\EK(\C t\oplus \g_{\pm }\oplus \C t^{*})\cong (\C[t]\ltimes \EK(\g_{\pm }))\otimes \C
[t^{*}]
\end{equation}
(which is a Hopf algebra isomorphism in the ``minus'' case), 
and define $c_{(-,\mp 1)}$ to be the one dimensional 
$\EK(\C t\oplus \g_{\pm }\oplus \C t^{*})$-module, 
where $\EK(\g_{\pm })$ and $t$ act trivially, and $t^{*}$ is defined to
act by $\mp 1$.   
Then let  $\Mb_{\pm}^{q}$ be the induced module  
$$
\Mb_{\pm }^{q}=\mathrm{Ind}_{\EK(\C t\oplus \g_\pm \oplus \C
t^*)}^{\EK(\g')}c_{(-,\mp 1)}.
$$
 
By construction, $\Mb_{\pm}^{q}$ is topologically free.  
Moreover, by the Poincare-Birkhoff-Witt Theorem we have 
an isomorphism of $k[[h]]$-modules
$\Mb_{\pm}^{q}=(U(\g_{\mp})[[h]])\oneq_{\pm}$ 
where $\oneq_{\pm}\in \Mb_{\pm}^{q}$.   

%  LEMMA 

\begin{lemma}\label{L:1}  We have the following isomorphisms of $\EK(\g')$-modules:
$F(M_{-})\cong M_{-}^{q}$, $F(M_+^*)=(M_+^q)^*$, 
$F(\Mb_{-})\cong \Mb_{-}^{q}$, and $F(\Mb_{+}^{*})\cong (\Mb_{+}^{q})^{*}$.
\end{lemma}

\begin{proof}
Let us prove that $F(M_{-})\cong M_{-}^{q}$; the other statements
are proved analogously.  Since $F(M_{-})$ and $M_{-}^{q}$ are
both deformations of $M_{-}$ as $U(\g')$-modules, it suffices to
show that $F(M_{-})$ contains a vector $v$ such that for all $a\in \EK(\g'_{-})$, 
$$av=\epsilon(a)v.$$
The action of $\EK(\g')$ on $F(M_{-})$ can be described as follows.  $$\EK(\g')=\EK(\g'_{+})\cdot\EK(\g'_{-}).$$  By definition $\EK(\g'_{+})\cong F(M_{-})$ acts by multiplication, and $\EK(\g'_{-})\cong \EK(\g'_{+})^{*cop}$ acts by comultiplication.  So the condition $av=\epsilon(a)v$ is equivalent to saying that $\Delta v = v\otimes 1$.  The element $v=1$ (unit of $F(M_{-})=\EK(\g'_{+})$) satisfies this condition.   
\end{proof}

We now define a representation of $\Pi_{n}$ coming from the
braiding of $\Mp$.   We must be careful here because the
associator of $\Mp$ is nontrivial.  Let  $\Mb_{+}^{*}\hat\otimes
\Umod^{\otimes n}$ be the object of $\Mp$ equipped with the
unique system of parentheses opening only at the extreme left (to
simplify notation we do not write the parentheses).  For  example
if $n=3$ we have $\Mb_{+}^{*}\hat\otimes \Umod^{\otimes 3}$ is
$((\Mb_{+}^{*}\hat\otimes \Umod)\otimes \Umod)\otimes \Umod$.  Let
$\rho_{1}$ be the representation of $\Pi_{n}$
on the space $Y:=\Hom_{\g_+\oplus \g_-\oplus \mathbb C t^*}(\Mb_-,\Mb_{+}^{*}\hat\otimes
\Umod^{\otimes n})$ defined by the formulas 
\begin{align*}
\label{}
   X_{0}& \rightarrow \big( \beta_{\Umod, \Mb_{+}^{*}}\otimes id^{\otimes (n-1)}\big)\circ\big(\beta_{\Mb_{+}^{*}, \Umod}\otimes id^{\otimes (n-1)}\big)  \\
  b_{i}  &  \rightarrow \Phi_{i}^{-1} \circ \beta_{i,i+1} \circ \Phi_{i} \;\text{ for } \; i=1,...,n-1
\end{align*}
where $\Phi_{i}$ is the morphism given by the action of the invertible element
$$\Phi_{i}=\dpi^{(i+2)}(\Phi)\otimes 1^{\otimes(n-i-1)}$$ 
expressed in terms of the map $\dpi^{(i+2)}:\Umod[[h]]^{\otimes 3}\rightarrow \Umod[[h]]^{\otimes (i+2)}$ defined recursively by $\dpi^{(3)}=id_{\Umod[[h]]^{\otimes 3}}$ and $\dpi^{(i+2)}=(\Delta\otimes id_{\Umod[[h]]}^{\otimes (i)})\dpi^{(i+1)}$.

Let $Y^q$ be the superspace 
\begin{equation}
\label{E:3}
Y^q=\Hom_{\EK(\g_{+}\oplus \g_{-}\oplus \C t^{*})}\left(\Mb_{-}^{q},(\Mb_{+}^{q})^{*}\hat\otimes \Umod^{\otimes n}[[h]]\right).
\end{equation}
Define $\rho_{2}$ to be the representation of $\Pi_{n}$ on $Y^q$ 
coming from the action of the R-matrix of $\EK(\g')$ on $(\Mb_{+}^{q})^{*}\hat\otimes \Umod^{\otimes n}[[h]]$. 

% PROP $\rho_{1}$ and $\rho_{2}$ are isomorphic 

\begin{proposition}\label{P:2}
The representations $\rho_{1}$ and $\rho_{2}$ are isomorphic.
\end{proposition}

\begin{proof}
  We have the following isomorphisms which preserve the braiding:
$$
\Hom_{\g_+\oplus \g_-\oplus \mathbb C t^*}(\Mb_{-},\Mb_{+}^{*}\hat\otimes \Umod^{\otimes n})[[h]]
$$
$$
\cong \Hom_{\EK(\g_+\oplus \g_-\oplus \mathbb C t^*)}
(F(\Mb_-), F( \Mb_{+}^{*}\hat\otimes
\Umod^{\otimes n})) 
$$
$$
\cong \Hom_{\EK(\g_+\oplus \g_-\oplus \mathbb C t^*)}(\Mb_{-}^{q},
(\Mb_{+}^q)^{*}\hat\otimes \Umod^{\otimes n}),
$$
where the first isomorphism comes from the fact 
that $F$ is fully faithful, and 
the second one from Lemma \ref{L:1}. This implies the
result. 

\end{proof}

% REMARK

Since the associator of $\Mp$ arises from the KZ equations, 
Proposition \ref{P:1} implies that the representation $\rho_1$
is the monodromy representation $\theta$ of the trigonometric KZ equations.  
Therefore, Proposition \ref{P:2} tells us that we can use the 
representation $\rho_{2}$ to compute $\theta$.  
We will now do this explicitly.  

Identifying the space 
$Y^q$
with $ \Umod^{\otimes n}[[h]]$ by Frobenius reciprocity, 
we see that the elements $b_i$ indeed act as
$\sigma_{ii+1}R_{ii+1}$. Thus it is enough for us to 
compute the action of $X_1$ upon this identification.  

We have $X_1=b_0^2$ in $B_{n+1}$.  Let $\Psi\in Y^q$.
Then we have $X_1\Psi=\tilde R_{10}\tilde R_{01}\Psi$, 
where $\tilde R$ is the universal R-matrix of $\EK(\g')$.  
In particular, $X_1\Psi \oneq_-=\tilde R_{10}\tilde R_{01}\Psi
\oneq_-$.  
To compute the value of $\tilde R_{10}\tilde R_{01}\Psi \oneq_-$ 
we need the following lemma and corollary.  

Recall that from functoriality of quantization we have 
a natural epimorphism of Hopf superalgebras 
$$\chi : \EK(\g_{+}') \rightarrow \EK(\C t^{*})=\C[t^{*}][[h]]$$.
 Let $\tilde R=\sum a_{i}\otimes b_{i}$. 

% LEMMA chi=e

\begin{lemma}\label{L:3}
$\sum \chi(a_{i})|_{t^{*}=-1}\otimes b_{i}=e^{-h(1\otimes t)}$.
\end{lemma}
\begin{proof}
To simplify notation set $A=\EK(\g'_{+})$ and $B=\EK(\C t^{*})$.
Let $A^{\star}$ be the dual QUE superalgebra of $A$, with
opposite coproduct.  As shown in
\cite{EK1,G04A}, there is a canonical isomorphism the QUE
superalgebras $A^{\star}$ and $\EK(\g'_{-})$, coming from 
the R-matrix.   We use this isomorphism to identify $A^{\star}$ and $\EK(\g'_{-})$ and write $A^{\star}=\EK(\g'_{-})$.  Similarly, $B^{\star}=\EK(\C t)$.

Let $\eta_{0}:\C t \rightarrow \g'_{-}=\g_{-}\oplus \C t$ be the obvious Lie superbialgebra embedding.  This morphism gives rise to a QUE superalgebra embedding 
\begin{equation}
\label{E:4}
\eta: \EK(\C t) \rightarrow \EK(\g'_{-})
\end{equation}
(the dual of $\chi$).
In the above notation we have $\eta: B^{\star}\rightarrow A^{\star}$.  Using this embedding we can regard $B^{\star}$ as a Hopf sub-superalgebra of $A^{\star}$.  

Let $D(A)=A\otimes A^{\star}$ be the quantum double of $A$.  As
shown in \cite{EK1,G04A}, $\EK(\g')$ is isomorphic to $D(A)$ as
braided QUE superalgebras.  In particular, $\tilde R$ corresponds
to the canonical element $R_A$ of $A\otimes A^*\subset 
D(A)\otimes D(A)$, and it suffices to show that
the statement of the lemma holds in $D(A)$, where $\tilde R$ is the canonical element.

Now, by the definition of the canonical element, 
$(\chi\otimes Id)(R_{A})$ belongs to $B\otimes B^{\star}$ and equals $R_{B}$, the universal R-matrix of $D(B)$.  But by definition $R_{B}=e^{ht^{*}\otimes t}$, and so $\sum \chi(a_{i})\otimes b_{i}=e^{ht^{*}\otimes t}$.  The lemma follows.
\end{proof}

%%  COROLLARY

\begin{corollary}\label{C:1}
$\sum \chi\left(S(a_{i})\right)|_{t^{*}=-1}\otimes b_{i}= e^{h(1\otimes t)}$.
\end{corollary}
\begin{proof}
Since $(S \otimes S)(\tilde R)=\tilde R$, we have 
$$\sum \chi\left(S(a_{1})\right)|_{t^{*}=-1}\otimes b_{i}=\sum (1\otimes S^{-1})(\chi\left(a_{1}\right)|_{t^{*}=-1}\otimes b_{i})=(1\otimes S^{-1})(e^{-h(1\otimes t)}).$$
This implies the corollary.
\end{proof}

%  LEMMA 

The main theorem now follows from the following lemma. 

\begin{lemma}\label{L:4}
 Set $T=(1\otimes S)(\tilde R)$.
 The representation $\rho_{2}$ is given by the formulas
\begin{align*}
\label{}
 b_{i}&\rightarrow \flip_{i,i+1}R_{i,i+1} &  X_{1} &\rightarrow e^{h(s+\rho_r)^{(1)}}m_{01}(T_{01}...T_{0n}) 
\end{align*}
 where $m$ is the multiplication in $\EK(\g')$.
\end{lemma}
\begin{proof}
The map $\phi:(\Mb_{+}^{q})^{*}\hat\otimes (\Umod[[h]])^{\hat\otimes
n} \rightarrow \Hom_{\C[[h]]}(\Mb_{+}^{q},(\Umod[[h]])^{\otimes
n})$ given by $\phi(f\otimes x)(y)=f(y)x$, $f\in
(\Mb_{+}^{q})^{*}$, $x\in (\Umod[[h]])^{\otimes n}$, $y\in
\Mb_{+}^{q}$, is an isomorphism.  Moreover, $\phi$ maps
$\EK(\g'_{-})$-invariant tensors to $\EK(\g'_{-})$-invariant
morphisms.   By Frobenius reciprocity, the space
$Hom_{\EK(\g'_{-})}(\Mb_{+}^{q},(\Umod[[h]])^{\otimes n})$ is
isomorphic to $(\Umod[[h]])^{\otimes n}$, via $g \mapsto
g(\oneq_{+})$.  Thus, the evaluation map $\big<,\big>:
((\Mb_{+}^{q})^{*}\hat\otimes (\Umod[[h]])^{\otimes n}
)^{\EK(\g'_{-})}\rightarrow (\Umod[[h]])^{\otimes n}$ determined
by $\big<\oneq_{+}\otimes 1^{\otimes n}, f\otimes x\big> =
f(\oneq_{+})x$, is an isomorphism.

Consider
 \begin{align}
\label{E:8}
  \big<\oneq_{+}\otimes 1^{\otimes n}, R_{10}R_{01}\Psi\oneq_{-}  \big>  & = \big<\oneq_{+}\otimes 1^{\otimes n}, \sum \big((-1)^{\p{a}_{i}\p{b}_{j}+\p{a}_{j}\p{b}_{j}}a_{i}b_{j}\otimes b_{i}a_{j}\otimes 1^{\otimes (n-1)}\big)\Psi\oneq_{-}\big> \notag \\
    &  = \sum \big< S(a_{i})\oneq_{+}\otimes 1^{\otimes n},  \big((-1)^{\p{a}_{i}\p{b}_{j}+\p{a}_{j}\p{b}_{j}}b_{j}\otimes b_{i}a_{j}\otimes 1^{\otimes (n-1)}\big)\Psi\oneq_{-}\big> \notag \\
    & = \sum \big<\oneq_{+}\otimes 1^{\otimes n},(-1)^{\p{a}_{j}\p{b}_{j}}b^{(0)}_{j}(e^{-ht}a_{j})^{(1)}\Psi\oneq_{-}\big>
\end{align}
 where the last equality comes from Corollary \ref{C:1}.  
 
Since $\Psi\oneq_{-}$ is $\EK(\g'_{-})$-invariant, as discussed above we can consider it as an element 
of $\Hom_{\EK(\g'_{-})}(\Mb_{+}^{q}, (\Umod[[h]])^{\otimes n})$.   In particular, we can identify $\sum(-1)^{\p{a}_{j}\p{b}_{j}}b^{(0)}_{j}(e^{-ht}a_{j})^{(1)}\Psi\oneq_{-}$ with
\begin{equation}
\label{E:5}
\sum (e^{ht}a_{j})^{(1)}\phi(\Psi\oneq_{-})\cdot S(b_{j})
\end{equation}
where $\cdot S(b_{j})$ means that $S(b_{j})$ first acts on $\Mb_{+}^{q}$.  Since $\phi(\Psi\oneq_{-})$ is $\EK(\g'_{-})$-invariant, we have that the element in \eqref{E:5} is equal to 
\begin{align}
\label{E:6}
  \sum
(e^{ht}a_{j})^{(1)}\Delta_{n}(S(b_{j}))\phi(\Psi\oneq_{-})  & =
(e^{ht})^{(1)}[(1\otimes S)(\tilde R)]^{1,12..n}\phi(\Psi\oneq_{-}) \notag \\
    &  = (e^{ht})^{(1)}T^{1,12..n}\phi(\Psi\oneq_{-}).
\end{align}
Here for $X\in \EK(\g)\otimes \EK(\g)$ we define $X^{1,12...n}$ to be 
$m_{12}((1\otimes \Delta_n)(X))$, where 
$m$ is the multiplication, and $\Delta_n: \EK(\g)\to \EK(\g)^{\otimes n}$ 
the iterated comultiplication. 
Using the hexagon relations for $\tilde R$, we can rewrite~(\ref{E:6}) as
\begin{align}
\label{E:7}
  \sum (e^{ht}a_{j})^{(1)}\Delta_{n}(S(b_{j}))\phi(\Psi\oneq_{-})  &=
(e^{ht})^{(1)}m_{01}(T_{01}...T_{0n})\phi(\Psi \oneq_-).
\end{align}
The lemma follows from combining above results.
\end{proof}

\end{document}